\noindent\centerline{\bf Erd\'elyi-Kober Fractional Integral Operators from a Statistical Perspective -I}

\vskip.3cm\centerline{A.M. MATHAI}
\vskip.2cm\centerline{Centre for Mathematical Sciences,}
\vskip.1cm\centerline{Arunapuram P.O., Pala, Kerala-68674, India, and}
\vskip.1cm\centerline{Department of Mathematics and Statistics, McGill University,}
\vskip.1cm\centerline{Montreal, Quebec, Canada, H3A 2K6}
\vskip.2cm\centerline{and}
\vskip.2cm\centerline{H.J. HAUBOLD}
\vskip.1cm\centerline{Office for Outer Space Affairs, United Nations}
\vskip.1cm\centerline{P.O. Box 500, Vienna International Centre} 
\vskip.1cm\centerline{A - 1400 Vienna, Austria, and}
\vskip.2cm\centerline{Centre for Mathematical Sciences,}
\vskip.1cm\centerline{Arunapuram P.O., Pala, Kerala-68674, India}

\vskip.5cm\noindent{\bf Abstract}

\vskip.3cm In this article we examine the densities of a product and a ratio of two real positive scalar random variables $x_1$ and $x_2$, which are statistically independently distributed, and we consider the density of the product $u_1=x_1x_2$ as well as the density of the ratio $u_2={{x_2}\over{x_1}}$ and show that Kober operator of the second kind is available as the density of $u_1$ and Kober operator of the first kind is available as the density of $u_2$ when $x_1$ has a type-1 beta density and $x_2$ has an arbitrary density. We also give interpretations of Kober operators of the second and first kind as Mellin convolution for a product and ratio respectively. Then we look at various types of generalizations of the idea thereby obtaining a large collection of operators which can all be called generalized Kober operators. One of the generalizations considered is the pathway idea where one can move from one family of operators to another family and yet another family and eventually end up with an exponential form. Common generalizations in terms of a Gauss' hypergeometric series is also given a statistical interpretation and put on a more general structure so that the standard generalizations given by various authors, including Saigo operators, are given statistical interpretations and are derivable as special cases of the general structure considered in this article.

\vskip.3cm\noindent{\bf 1.\hskip.3cm Introduction}

\vskip.3cm Let $f_1(x_1)$ be the density of $x_1$ and $f_2(x_2)=f(x_2)$ be the density of $x_2$ where let $x_1$ and $x_2$ be independently distributed with $x_1$ having the density
$$f_1(x_1)={{\Gamma(\zeta+1+\alpha)}\over{\Gamma(\zeta+1)\Gamma(\alpha)}}x_1^{\zeta}(1-x_1)^{\alpha-1},~0<x_1<1,\Re(\alpha)>0,\Re(\zeta)>-1
$$and $f_1(x_1)=0$ elsewhere. In statistical problems, usually the parameters are real but the results will hold for complex parameters and hence we list the conditions for complex parameters. When $x_1$ and $x_2$ are independently distributed, from standard procedure of transformation of variables, the density of $u=x_1x_2$ can be written as
$$g(u)=\int_v{{1}\over{v}}f_1({{u}\over{v}})f_2(v){\rm d}v=\int_y{{1}\over{y}}f_1(y)f_2({{u}\over{y}}){\rm d}y.\eqno(1.1)
$$Taking $f_1(x_1)$ as type-1 beta and $f_2(x_2)=f(x_2)$ as arbitrary
$$\eqalignno{g(u)&={{\Gamma(\alpha+\zeta+1)}\over{\Gamma(\zeta+1)\Gamma(\alpha)}}\int_{t=u}^{\infty}{{1}\over{t}}({{u}\over{t}})^{\zeta}(1-{{u}\over{t}})^{\alpha-1}f(t){\rm d}t\cr
&={{\Gamma(\alpha+\zeta+1)}\over{\Gamma(\zeta+1)}}{{u^{\zeta}}\over{\Gamma(\alpha)}}\int_{u}^{\infty}(t-u)^{\alpha-1}t^{-\zeta-\alpha}f(t){\rm d}t\cr
&={{\Gamma(\alpha+\zeta+1)}\over{\Gamma(\zeta+1)}}K_{u}^{\zeta,\alpha}f.\cr
\noalign{\hbox{This means that}}
K_{u}^{\zeta,\alpha}f&={{\Gamma(\zeta+1)}\over{\Gamma(\alpha+\zeta+1)}}g(u).&(1.2)\cr}
$$Thus, we have the following theorem:

\vskip.3cm\noindent{\bf Theorem 1.1.}\hskip.3cm{\it Kober fractional integral operator of the second kind is a constant multiple of the density of a product of two real scalar statistically independently distributed positive random variables $x_1$ and $x_2$ where $x_1$ has a type-1 beta density with the parameters $(\zeta+1,\alpha)$, and $x_2$ has an arbitrary density $f(x_2)$.}

\vskip.3cm Since $x_1$ and $x_2$ are independently distributed, we have
$$E(u^{s-1})=E(x_1^{s-1})E(x_2^{s-1})
$$where $E(\cdot)$ denotes the expected value. But
$$E(x_1^{s-1})={{\Gamma(\alpha+\zeta+1)}\over{\Gamma(\zeta+1)}}{{\Gamma(\zeta+s)}\over{\Gamma(\alpha+\zeta+s)}}\hbox{  for  }\Re(\zeta)>0,\Re(\zeta+s)>0,\Re(\alpha)>0
$$and let $E(x_2^{s-1})=f^{*}(s)= $ the Mellin transform of $f(x)$. If the Mellin transform of the $g(u)$, with Mellin parameter $s$, is denoted by $M\{g(u);s\}$ then
$$M\{{{\Gamma(\zeta+1)}\over{\Gamma(\alpha+\zeta+1)}}g(u);s\}={{\Gamma(\zeta+s)}\over{\Gamma(\alpha+\zeta+s)}}f^{*}(s)=M\{K_{u}^{\zeta,\alpha}f;s\}.\eqno(1.3)
$$Hence Kober fractional integral operator of the second kind can be considered as a Mellin convolution for a product. Then, naturally the inverse Mellin transform of (1.3) provides explicit expression for Kober fractional integral operator of the second kind, namely,
$$K_{u}^{\zeta,\alpha}f={{1}\over{2\pi i}}\int_{c-i\infty}^{c+i\infty}{{\Gamma(\zeta+s)}\over{\Gamma(\alpha+\zeta+s)}}f^{*}(s)u^{-s}{\rm d}x\eqno(1.4)
$$where the form is available through the convolution integral coming from the Mellin convolution of a product, that is, a type-1 beta form convoluted with the arbitrary function $f(x)$.

\vskip.3cm\noindent{\bf 1.1.\hskip.3cm A Pathway Kober Operator of the Second Kind}

\vskip.3cm Let $f_1(x_1)$ be the pathway density
$$\eqalignno{f_1(x_1)&=c_1~x_1^{\gamma}[1-a(1-q)x_1^{\delta}]^{{\eta}\over{1-q}}&(1.5)\cr
\noalign{\hbox{for $q<1,\eta>0,a>0,\delta>0$ where}}
c_1&={{\delta[a(1-q)]^{{\gamma+1}\over{\delta}}\Gamma({{\gamma+1}\over{\delta}}+{{\eta}\over{1-q}}+1)}
\over{\Gamma({{\gamma+1}\over{\delta}})\Gamma({{\eta}\over{1-q}}+1)}}.\cr}
$$Then the density of $u=x_1x_2$, where $x_2$ has arbitrary density $f(x_2)$, is given by
$$\eqalignno{g(u)&=c_1\int_v{{1}\over{v}}f_1({{u}\over{v}})f_2(v){\rm d}v\cr
&=c_1\int_v{{1}\over{v}}({{u}\over{v}})^{\gamma}[1-a(1-q)({{u}\over{v}})^{\delta}]^{{\eta}\over{1-q}}f(v){\rm d}v&(1.6)\cr
&=c_1u^{\gamma}\int_{v=u[a(1-q)]^{{1}\over{\delta}}}^{\infty}\{[v^{\delta}-a(1-q)u^{\delta}]^{{\eta}\over{1-q}}
 v^{-\gamma-({{\eta\delta}\over{1-q}}+1)}\}f(v){\rm d}v&(1.7)\cr}
$$
\vskip.3cm\noindent{\bf 1.2.\hskip.3cm Special Cases}
\vskip.3cm\noindent{\bf Case (1):}\hskip.3cm When $\delta =m,q=0,{{\eta}\over{1-q}}=\alpha-1$ the right side of (1.7) agrees with the result (2.6.9) of Mathai and Haubold (2008).

\vskip.3cm\noindent{\bf Case (2):}\hskip.3cm When $\delta=1,a=1,q=0,\eta=\alpha-1$ then we have
$$\eqalignno{{{\Gamma(\gamma+1)}\over{\Gamma(\gamma+1+\alpha)}}g(u)&={{1}\over{\Gamma(\alpha)}}u^{\gamma}
\int_{v=u}^{\infty}(v-u)^{\alpha-1}v^{-\gamma-\alpha}f(v){\rm d}v\cr
&=K_{u}^{\gamma,\alpha}f=\hbox{ Kober operator of the second kind}&(1.8)\cr}
$$
\vskip.3cm\noindent{\bf Case (3):}\hskip.3cm When $\delta=1,a=1,q=0,\eta=\alpha-1,\gamma=0$ then
$$\eqalignno{{{1}\over{\Gamma(\alpha+1)}}g(u)&={{1}\over{\Gamma(\alpha)}}\int_{t=u}^{\infty}(t-u)^{\alpha-1}
t^{-\alpha}f(t){\rm d}t\cr
&=K_{u}^{0,\alpha}f=\hbox{ special case of Kober operator}&(1.9)\cr
&={_xW_{\infty}^{-\alpha}}t^{-\alpha}f(t)&(1.10)\cr}
$$is the Weyl right sided fractional integral operator of order $\alpha$ or right sided Riemann-Liouville fractional integral operator of order $\alpha$, when the right limit is at $\infty$, for the function $t^{-\alpha}f(t)$.
\vskip.3cm When $q$ moves from $-\infty$ to $1$ then (1.7) describes a collection of generalized Kober operators of the second kind, operating on an arbitrary function $f(t)$. It can also be considered as a Mellin convolution of a product where one function $f(x_2)$ is arbitrary and the other function $f_1(x_1)$ is of the form in (1.5). Here $q$ describes a path of movement of the Kober operator of the second kind. In the limit when $q\to 1_{-}$ then (1.7) will go to

$$\eqalignno{\lim_{q\to 1_{-}}g(u)&=c_1^{*}\int_{v=0}^{\infty}{{1}\over{v}}({{u}\over{v}})^{\gamma}{\rm e}^{-a\eta({{u}\over{v}})^{\delta}}f(v){\rm d}v&(1.11)\cr
&=c_1^{*}u^{\gamma}\int_{v=0}^{\infty}v^{-\gamma-1}{\rm e}^{-a\eta({{u}\over{v}})^{\delta}}f(v){\rm d}v\cr
\noalign{\hbox{where}}
c_1^{*}&=\delta{{(a\eta)^{{\gamma+1}\over{\delta}}}\over{\Gamma({{\gamma+1}\over{\delta}})}}.\cr}
$$The pathway form of (1.11) is also connected to Kr\"atzel transform if $f(v)$ can be written as ${\rm e}^{-bv}\phi(v)$. Then the integral in (1.11) will correspond to generalized Kr\"atzel transform of $\phi(v)$. There are lots of applications of Kr\"tzel transform in various disciplines. This transform is also connected to inverse Gaussian density in stochastic processes, to Bayesian analysis, reaction rate probability integral in reaction rate theory and many other topics, the details may be seen from Mathai (2012), Kumar (2010), Kumar and Kilbas (2010) and Kumar and Haubold (2010).

When $q>1$ then writing $1-q=-(q-1)$ with $q>1$, $f_1(x_1)$ of (1.5) changes to the following form. For $q>1$,
$$\eqalignno{f_1(x_1)&=c_2x_1^{\gamma}[1+a(q-1)x_1^{\delta}]^{-{{\eta}\over{q-1}}},a>0,\delta>0,\eta>0,q>1&(1.12)\cr
\noalign{\hbox{where}}
c_2&=\delta{{[a(q-1)]^{{\gamma+1}\over{\delta}}\Gamma({{\eta}\over{q-1}})}\over{\Gamma({{\gamma+1}\over{\delta}})
\Gamma({{\eta}\over{q-1}}-{{\gamma+1}\over{\delta}})}},~{{\eta}\over{q-1}}-{{\gamma+1}\over{\delta}}>0&(1.13)\cr
&\to \delta{{(a\eta)^{{\gamma+1}\over{\delta}}}\over{\Gamma({{\gamma+1}\over{\delta}})}}\hbox{  when }q\to 1_{+}.\cr}
$$Note that in this case $0<x_1<\infty$. Then proceeding as before

$$\eqalignno{g(u)&=c_2\int_{t=0}^{\infty}{{1}\over{t}}({{u}\over{t}})^{\gamma}[1+a(q-1)({{u}\over{t}})^{\delta}]^{-{{\eta}\over{q-1}}}f(t){\rm d}t\cr
&=c_2u^{\gamma}\int_0^{\infty}t^{-\gamma+({{\delta\eta}\over{q-1}}-1)}[t^{\delta}+a(q-1)u^{\delta}]^{-{{\eta}\over{q-1}}}f(t){\rm d}t.&(1.14)\cr}
$$This can also be considered as a generalization of the Kober operator of the second kind. It goes into the following form when $q\to1_{+}$.
$$\lim_{q\to 1_{+}}g(u)=c_2^{*}\int_0^{\infty}t^{-\gamma-1}{\rm e}^{-a\eta({{u}\over{t}})^{\delta}}f(t){\rm d}t\eqno(1.15)
$$This is the same as (1.11). We can show that
$$\lim_{q\to 1_{-}}c_1^{*}=\lim_{q\to 1_{+}}c_2^{*}.
$$
\vskip.3cm\noindent{\bf 1.3.\hskip.3cm Another Form of Generalization of Kober Operators of the Second Kind}

\vskip.3cm Instead of taking the type-1 beta density for $f_1(x_1)$ and an arbitrary density for $f_2(x_2)$, one can take any given density for $f_1(x_1)$ and an arbitrary density for $x_2$ and then take the Mellin convolution of a product. Then this will give a class of generalized Kober operators from a statistical point of view. If a fractional type integral is required then the variable can be relocated at $x=b$ so that $x\ge b$ for some $b$. From  a mathematical point of view such a generalization may not have much of a significance.

\vskip.3cm\noindent{\bf 1.4.\hskip.3cm A Generalization in terms of Hypergeometric Series}

\vskip.3cm Let us consider appending a hypergeometric series to  the basic density of $x_1$. Consider the function $$f_1(x_1)={{1}\over{c}}~{_pF_q}(a_1,...,a_p;b_1,...,b_q ;ax_1)x_1^{\zeta}(1-x_1)^{\alpha-1},~0<x_1<1\eqno(1.16)
$$and $f_1(x_1)=0$ elsewhere, where $c$ is the normalizing constant. We can create a statistical density out of this form as follows: In order to assure nonnegativity of the function let us assume that the parameters $a_1,...,a_p,b_1,...,b_q,a$ are all  positive. Then
$$\eqalignno{{_pF_q}(a_1,...,a_p;&b_1,...,b_q:ax_1)x_1^{\zeta}(1-x_1)^{\alpha-1}\cr
&=\sum_{k=0}^{\infty}{{(a_1)_k...(a_p)_k}\over{(b_1)_k...(b_q)_k}}{{a^kx_1^k}\over{k!}}x_1^{\zeta}(1-x_1)^{\alpha-1}.\cr}
$$Total integral is available from the basic integral
$$\int_0^1x_1^{\zeta+k}(1-x_1)^{\alpha-1}{\rm d}x_1={{\Gamma(\zeta+1+k)\Gamma(\alpha)}\over{\Gamma(\alpha+\zeta+1+k)}}={{\Gamma(\alpha)\Gamma(\zeta+1)}\over{\Gamma(\alpha+\zeta+1)}}{{(\zeta_1)_k}\over{(\alpha+\zeta+1)_k}}.
$$Then the normalizing constant
$$c={{\Gamma(\alpha)\Gamma(\zeta+1)}\over{\Gamma(\alpha+\zeta+1)}}{_{p+1}F_{q+1}}(a_1,...,a_p,\zeta+1;
b_1,....,b_q,\alpha+\zeta+1;a).
$$Then
$$f_1(x_1)={{1}\over{c}}~{_pF_q}(a_1,...,a_p;b_1,...,b_q;ax_1)x_1^{\zeta}(1-x_1)^{\alpha-1},0<x_1<1
$$and zero elsewhere is a density. Take this form of $f_1(x_1)$ and proceed to find the density of $u=x_1x_2$ as before. Denoting the density, again by $g(u)$,
$$\eqalignno{g(u)&={{1}\over{c}}\sum_{k=0}^{\infty}{{(a_1)_k...(a_p)_k}\over{(b_1)_k...(b_q)_k}}{{a^k}\over{k!}}
\int_vu^{\zeta}(v-u)^{\alpha-1}v^{-\zeta-\alpha}({{u}\over{v}})^kf(v){\rm d}v\cr
&={{1}\over{c}}u^{\zeta}\int_{v>u}(v-u)^{\alpha-1}v^{-\zeta-\alpha}
{_pF_q}(a_1,...,a_p;b_1,...,b_q;{{au}\over{v}})f(v){\rm d}v.&(1.17)\cr}
$$A particular case of this for a ${_2F_1}$ is equation (2.7.2) of Mathai and Haubold (2008). This particular case was given by others earlier. Note that there is one main drawback in taking a ${_2F_1}$ because then there may be problems in taking Laplace, Mellin and other transforms for the convergence of the series forms. Hence it is safer to take $q\ge p$ in the case of appending a hypergeometric series to the type-1 beta form for $f_1(x_1)$. Note that (1.17) is a generalization of Kober operator of the second kind as well as one has an interpretation in terms of a statistical density.
\vskip.2cm Another form of appending a hypergeometric series is to consider a hypergeometric series with argument $a(1-x_1)$ instead $ax_1$. Going through the same process as before, one can create a statistical density of the form
$$f_1(x_1)={{1}\over{\tilde{c}}}\sum_{k=0}^{\infty}{{(a_1)_k...(a_p)_k}\over{(b_1)_k...(b_q)_k}}{{a^k}\over{k!}}x_1^{\zeta}(1-x_1)^{\alpha-1+k}\eqno(1.18)
$$for $0<x_1<1$ and zero elsewhere, where
$$\tilde{c}={{\Gamma(\zeta+1)\Gamma(\alpha)}\over{\Gamma(\alpha+\zeta+1)}}{_{p+1}F_{q+1}}(a_1,...,a_p,\alpha;
b_1,...,b_q,\zeta+1+\alpha;a).
$$In order to guarantee nonnegativity we may assume all parameters $a_j$'s, $b_j$'s, be positive, $a>0$, $\alpha>0,q\ge p$. If $p=q+1$ then take $|a(1-x_1)|<1$. Proceeding exactly as before, taking $x_1$ having this appended density and $x_2$ having an arbitrary density, then the density of $u=x_1x_2$, again denoted by $g(u)$, is available as
$$\eqalignno{g(u)&={{1}\over{\tilde{c}}}\sum_{k=0}^{\infty}{{(a_1)_k...(a_p)_k}\over{(b_1)_k...(b_q)_k}}{{a^k}\over{k!}}
 \int_v{{1}\over{v}}({{u}\over{v}})^{\zeta}(1-{{u}\over{v}})^{\alpha+k-1}f(v){\rm d}v.\cr}
$$The integral part reduces to
$$u^{\zeta}\int_{v>u}(v-u)^{\alpha-1}v^{-\zeta-\alpha}(1-{{u}\over{v}})^kf(v){\rm d}v.
$$Hence
$$\eqalignno{g(u)&={{u^{\zeta}}\over{\tilde{c}}}\int_{v>u}(v-u)^{\alpha-1}v^{-\zeta-\alpha}
{_pF_q}(a_1,...,a_p;b_1,...,b_q;a(1-{{u}\over{v}}))f(v){\rm d}v,~|a|<1,v>u.&(1.19)\cr}
$$This is a generalization of Kober operator of the second kind. Here also, one could have taken the argument of ${_pF_q}$ as $a^{\delta_1}(1-x_1)^{\delta_2}$. These will provide  more generalized forms. A particular case of (1.19) is (2.7.4) of Mathai and Haubold (2008). This particular case, in terms of a ${_2F_1}$, was given by others earlier. As remarked above, there is a disadvantage in taking a ${_2F_1}$. This special case in terms of a ${_2F_1}$ is Saigo operator, see (2.7.8) of Mathai and Haubold (2008).

\vskip.3cm\noindent{\bf Remark 1.1.}\hskip.3cm From the procedures in (1.17) and (1.19) it is clear that one can consider $f_1(x_1)$ in terms of a hypergeometric function with argument $ax_1$ or $a^{\delta_1}x_1^{\delta_2}$ or $a(1-x_1)$ or $a^{\delta_1}(1-x_1)^{\delta_2}$ or $a^{\delta_1}(1-x_1)^{\delta_2}x_1^{\delta_3}$ with $\delta_j>0,j=1,2,3$. The procedure will be the same. If statistical densities are not needed then one can take  multiplicative factors for $f_1(x_1)$ as well as for $f_2(x_2)$. Instead of a hypergeometric series, one can consider $f_1(x_1)$ in terms of a Meijer's G-function or H-function with arguments any one of them mentioned above. If $g(u)$ to remain as a statistical density then, apart from convergence of the series and integrals, the parameters are to be restricted so that the functions remain positive in the range $0<x_1<1$ and zero outside this range. Since these generalizations are routine mathematical exercises we will not give the explicit expressions for each generalization of Kober operator of the second kind here.

\vskip.3cm\noindent{\bf 1.5.\hskip.3cm Mellin Transform of the Generalized Kober Operator of the Second Kind}

\vskip.3cm For the generalized form in (1.17) the Mellin transform is available by evaluating the integral
$$\eqalignno{\int_0^{\infty}u^{s-1}u^{\zeta+k}&[\int_{v>u}(v-u)^{\alpha-1}v^{-\zeta-\alpha-k}f(v){\rm d}v]\cr
&=\int_{v=0}^{\infty}v^{-\zeta-\alpha-k}f(v)[\int_{u=0}^vu^{s-1+\zeta+k}(v-u)^{\alpha-1}{\rm d}u]{\rm d}v\cr
&=\int_{v=0}^{\infty}v^{s-1}f(v){\rm d}v\int_0^1y^{s+\zeta+k-1}(1-y)^{\alpha-1}{\rm d}y\cr
&={{\Gamma(\alpha)\Gamma(\zeta+s)}\over{\Gamma(\alpha+\zeta+s)}}{{(\zeta+s)_k}\over{(\alpha+\zeta+s)_k}}f^{*}(s).\cr}
$$Therefore the Mellin transform of (1.17) is the following:
$$\eqalignno{M\{g(u)\hbox{ of (1.17)};s\}&={{\Gamma(\alpha)}\over{c}}{{\Gamma(\zeta+s)}\over{\Gamma(\alpha+\zeta+s)}}\cr
&\times{_{p+1}F_{q+1}}(a_1,...,a_p,\zeta+s;b_1,...,b_q,\alpha+\zeta+s;a).&(1.20)\cr}
$$In a similar manner one can compute the Mellin transform of $g(u)$ of (1.19). The base integral to be evaluated is
$$\eqalignno{\int_{u=0}^{\infty}u^{s-1+\zeta}&[\int_{v=u}^{\infty}(v-u)^{\alpha-1}v^{-\zeta-\alpha}(1-{{u}\over{v}})^k{\rm d}u]f(v){\rm d}v\cr
&=\int_{v=0}^{\infty}v^{-\zeta-\alpha+\alpha-1}f(v)\cr
&\times[\int_{u=0}^v(1-{{u}\over{v}})^{\alpha+k-1}u^{\zeta+s-1}{\rm d}u]{\rm d}v\cr
&=\int_{v=0}^{\infty}v^{s-1}f(v){\rm d}v[\int_0^1y^{\zeta+s-1}(1-y)^{\alpha+k-1}{\rm d}y]\cr
&=f^{*}(s){{\Gamma(\zeta+s)\Gamma(\alpha)}\over{\Gamma(\alpha+\zeta+s)}}{{(\alpha)_k}\over{(\alpha+\zeta+s)_k}},\Re(\alpha)>0,\Re(\zeta+s)>0.\cr}
$$Therefore
$$\eqalignno{M\{g(u)\hbox{ of (1.19)};s\}&={{\Gamma(\alpha)}\over{\tilde{c}}}{{\Gamma(\zeta+s)}\over{\Gamma(\alpha+\zeta+s)}}\cr
&\times{_{p+1}F_{q+1}}(a_1,...,a_p,\alpha;b-1,...,b_q,\alpha+\zeta+s;a).&(1.21)\cr}
$$
\vskip.3cm\noindent{\bf 2.\hskip.3cm Kober Operator of the First Kind}

\vskip.3cm Let $x_1$ and $x_2$ be statistically independently distributed real positive scalar random variables. Let $u={{x_2}\over{x_1}}$. Let $x_1$ have a type-1 beta density with parameters $(\zeta,\alpha)$, that is, the density of $x_1$, denoted by $f_1(x_1)$, is given by
$$f_1(x_1)={{\Gamma(\zeta+\alpha)}\over{\Gamma(\zeta)\Gamma(\alpha)}}x_1^{\zeta-1}(1-x_1)^{\alpha-1},~0<x_1<1,
\Re(\alpha)>0,\Re(\zeta)>0.\eqno(2.1)
$$Let $x_2$ have an arbitrary density $f_2(x_2)=f(x_2)$ for some density $f(x_2)$. Then the density of $u={{x_2}\over{x_1}}$ is available by considering the transformation $u={{x_2}\over{x_1}},v=x_2$. Then ${\rm d}x_1\wedge{\rm d}x_2=-{{v}\over{u^2}}{\rm d}u\wedge{\rm d}v$. The joint density of $u$ and $v$ and from there the marginal density of $u$, again denoted as $g(u)$, is available as
$$g(u)=\int_vf_1({{v}\over{u}})f_2(v)(-{{v}\over{u^2}}){\rm d}v.\eqno(2.2)
$$Limits of $u$ will be from $\infty$ to $v$ and $0<v<u$. Then the marginal density is available as
$$\eqalignno{g(u)&=\int_{v=0}^uf_1({{v}\over{u}})f(v){{v}\over{u^2}}{\rm d}v\cr
&={{\Gamma(\zeta+\alpha)}\over{\Gamma(\zeta)\Gamma(\alpha)}}\int_{v=0}^u({{v}\over{u}})^{\zeta-1}
(1-{{v}\over{u}})^{\alpha-1}{{v}\over{u^2}}f(v){\rm d}v.\cr
\noalign{\hbox{Therefore}}
{{\Gamma(\zeta)}\over{\Gamma(\zeta+\alpha)}}g(u)&={{1}\over{\Gamma(\alpha)}}u^{-\zeta-\alpha}
\int_{v=0}^u(u-v)^{\alpha-1}v^{\zeta}f(v){\rm d}v&(2.3)\cr
&=I_{u}^{\zeta,\alpha}f.&(2.4)\cr}
$$
\vskip.3cm\noindent{\bf Theorem 2.1.}\hskip.3cm{\it The Kober operator of the first kind
$$I_{u}^{\zeta,\alpha}f={{u^{-\zeta-\alpha}}\over{\Gamma(\alpha)}}\int_{v=0}^u(u-v)^{\alpha-1}v^{\zeta}f(v){\rm d}v\eqno(2.5)
$$for $\Re(\zeta)>0,\Re(\alpha)>0$ is available as ${{\Gamma(\zeta)}\over{\Gamma(\zeta+\alpha)}}g(u)$ where $g(u)$ is the density of $u={{x_2}\over{x_1}}$ with $x_1$ having a type-1 beta density with parameters $(\zeta,\alpha)$ and $x_2$ has an arbitrary density, and $x_1$ and $x_2$ are statistically independently distributed.}

\vskip.3cm When $u={{x_2}\over{x_1}}$ where $x_1$ and $x_2$ are independently distributed, we have, denoting the expected value of $(\cdot)$ by $E(\cdot)$,
$$\eqalignno{E(u^{s-1})&=E(x_2^{s-1})E({{1}\over{x_1}})^{s-1}=E(x_2^{s-1})E(x_1^{-s+1})\cr
&=f^{*}(s){{\Gamma(\zeta+\alpha)}\over{\Gamma(\zeta)}}{{\Gamma(\zeta+1-s)}\over{\Gamma(\alpha+\zeta+1-s)}}&(2.6)\cr}
$$for $\Re(\alpha)>0,\Re(\zeta+1-s)>0$ where $f^{*}(s)$ is the Mellin transform of the arbitrary function $f(t)$. Therefore
$${{\Gamma(\zeta)}\over{\Gamma(\zeta+\alpha)}}g^{*}(x)=f^{*}(s)f_1^{*}(2-s).\eqno(2.7)
$$This means that the Kober operator of the first kind $I_{u}^{\zeta,\alpha}f$ can be considered as the Mellin convolution of a ratio $u={{x_2}\over{x_1}}$. Then Kober operator of the first kind can be considered as Mellin convolution of a ratio whereas that of the second kind as a Mellin convolution of a product.
\vskip.2cm Note that
$$\Gamma(\alpha)u^{\zeta+\alpha}I_{u}^{\zeta,\alpha}f=\int_{v=0}^uv^{\zeta}(u-v)^{\alpha-1}f(v){\rm d}v.
$$This is Euler transform of $f(v)$, see for example, Mathai et al.(2010). Note that all transforms where the basic function is a type-1 beta form or either of the form $x^{\alpha}(1-x)^{\beta},~0<x<1$ or of the form $(x-a)^{\alpha}(b-x)^{\beta},~a<x<b$ can be connected to Kober operators.

\vskip.3cm\noindent{\bf 2.1.\hskip.3cm A Pathway Generalization of Kober Operator of the First Kind}

\vskip.3cm Let the density of $x_1$ be of the form
$$\eqalignno{f_1(x_1)&=c_1x_1^{\gamma-1}[1-a(1-q)x_1^{\delta}]^{{\eta}\over{1-q}},1-a(1-q)x^{\delta}>0,\eta>0,q<1,\delta>0,a>0&(2.8)\cr
\noalign{\hbox{where}}
c_3&={{\delta[a(1-q)]^{{\gamma}\over{\delta}}\Gamma({{\eta}\over{1-q}}+1+{{\gamma}\over{\delta}})}\over{\Gamma({{\gamma}\over{\delta}})\Gamma({{\eta}\over{1-q}}+1)}}.\cr
\noalign{\hbox{Then}}
g(u)&=\int_vf_1({{v}\over{u}})f_2(v){{v}\over{u^2}}{\rm d}v\cr
&=c_3\int_v({{v}\over{u}})^{\gamma-1}[1-a(1-q)({{v}\over{u}})^{\delta}]^{{\eta}\over{1-q}}f_2(v){{v}\over{u^2}}{\rm d}v&(2.9)\cr
&=c_3u^{-\gamma-({{\delta\eta}\over{1-q}}+1)}
\int_vv^{\gamma}[u^{\delta}-a(1-q)v^{\delta}]^{{\eta}\over{1-q}}f(v){\rm d}v.&(2.10)\cr}
$$Through $q$ one has a collection of operators from (2.10), which can be treated as a generalization of Kober operator of the first kind.
\vskip.3cm\noindent{\bf 2.2.\hskip.3cm Some Special Cases}

\vskip.3cm\noindent{\bf Case (2.1):}\hskip.3cm For $a=1,q=0,\delta=1,{{\eta}\over{1-q}}=\alpha-1$ we have, calling this special case $g_1(u)$,
$$\eqalignno{g_1(u)&={{\Gamma(\gamma+\alpha)}\over{\Gamma(\gamma)\Gamma(\alpha)}}u^{-\gamma-\alpha}
\int_{v=0}^uv^{\gamma}(u-v)^{\alpha-1}f(v){\rm d}v&(2.11)\cr
\noalign{\hbox{for $\Re(\alpha)>0,\Re(\gamma)>0$. Therefore}}
{{\Gamma(\gamma)}\over{\Gamma(\alpha+\gamma)}}g_1(u)&={{u^{-\gamma-\alpha}}\over{\Gamma(\alpha)}}
\int_{v=0}^uv^{\gamma}(u-v)^{\alpha-1}f(v){\rm d}v\cr
&=I_{u}^{\gamma,\alpha}f.&(2.12)\cr}
$$Hence (2.10) gives a generalization of Kober operator of the first kind. This generalization also gives a path through $q$. For various values of $q$ one has a collection of functions which can all be considered as generalizations of Kober operators of the first kind. Note that (2.10) can also be looked upon as a generalization of the pathway integral transform introduced by Nair (2009). We can also look upon (2.10) as a generalized pathway transform of $f(x)$. Finally, when $q\to 1_{-}$ we have the following form:

\vskip.3cm\noindent{\bf Case (2.2):\hskip.3cm}
$$\lim_{q\to 1_{-}}g(u)=\delta{{(a\eta)^{{\gamma}\over{\delta}}}\over{\Gamma({{\gamma}\over{\delta}})}}u^{-\gamma-1}
\int_{v=0}^{\infty}v^{\gamma}{\rm e}^{-a\eta({{v}\over{u}})^{\delta}}f(v){\rm d}v.\eqno(2.13)
$$This integral part when $\delta=1$ is the Laplace transform of $v^{\gamma}f(v)$ with Laplace parameter ${{a\eta}\over{u}}$. \vskip.2cm All the functions described from (2.10) to (2.13) can be taken as generalizations of Kober operators of the first kind.

\vskip.3cm\noindent{\bf Case (2.3):}\hskip.3cm For $a=1,q=0,\delta=m,{{\eta}\over{1-q}}=\alpha-1$ we have a special case, call it $g_2(u)$. Then
${{\Gamma(\gamma)}\over{\Gamma(\gamma+\alpha)}}g_2(u)$ is (2.6.8) of Mathai and Haubold (2008).

\vskip.3cm\noindent{\bf 2.3.\hskip.3cm Generalization of Kober Operator of the First Kind in Terms of Hypergeometric Series}

\vskip.3cm Let us append a hypergeometric series to our basic function $x_1^{\zeta-1}(1-x_1)^{\alpha-1}$. Consider the function
$$\eqalignno{{_pF_q}(a_1,...,a_p&;b_1,...,b_q;ax_1)x_1^{\zeta-1}(1-x_1)^{\alpha-1}\cr
&=\sum_{k=0}^{\infty}{{(a_1)_k...(a_p)_k}\over{(b_1)_k...(b_q)_k}}{{a^kx_1^k}\over{k!}}x_1^{\zeta-1}(1-x_1)^{\alpha-1}\cr}
$$for $0<x_1<1$ and zero elsewhere. The basic integral part
$$\int_0^1x_1^{\zeta+k-1}(1-x_1)^{\alpha-1}{\rm d}x_1={{\Gamma(\zeta+k)\Gamma(\alpha)}\over{\Gamma(\alpha+\zeta+k)}}
={{\Gamma(\alpha)\Gamma(\zeta)}\over{\Gamma(\alpha+\zeta)}}{{(\zeta)_k}\over{(\alpha+\zeta)_k}}.\eqno(2.14)
$$Let
$$c^{(1)}={_{p+1}F_{q+1}}(a_1,...,a_p,\zeta;b_1,...,b_q,\alpha+\zeta;a){{\Gamma(\alpha)\Gamma(\zeta)}\over{\Gamma(\alpha+\zeta)}}.\eqno(2.15)
$$In order to guarantee positivity, let us assume that all parameters $a_j$'s, $b_j$'s and $a$ to be positive and $\alpha>0,\zeta>0$. Then the density of $u={{x_2}\over{x_1}}$, denoted again by $g(u)$, is given by
$$\eqalignno{g(u)&={{1}\over{c^{(1)}}}\sum_{k=0}^{\infty}{{(a_1)_k...(a_p)_k}\over{(b_1)_k...(b_q)_k}}{{a^k}\over{k!}}\cr
&\times\int_v({{v}\over{u}})^{\zeta+k-1}(1-{{v}\over{u}})^{\alpha-1}{{v}\over{u^2}}f(v){\rm d}v.\cr
\noalign{\hbox{The integral part is}}
u^{-\zeta-\alpha}&\int_{v=0}^uv^{\zeta}(u-v)^{\alpha-1}({{v}\over{u}})^kf(v){\rm d}v.\cr
g(u)&={{1}\over{c^{(1)}}}u^{-\zeta-\alpha}\int_{v=0}^uv^{\zeta}(u-v)^{\alpha-1}\cr
&\times{_pF_q}(a_1,...,a_p;b_1,...,b_q;a{{v}\over{u}})f(v){\rm d}v.&(2.16)\cr}
$$A particular case of (2.16) in terms of a ${_2F_1}$ is (2.7.1.) of Mathai and Haubold (2008). This particular case was given by others earlier as generalization of Kober operator of the first kind, not as a statistical density or as a Mellin convolution of a ratio. In the hypergeometric function we may replace the argument $ax_1$ by $a^{\delta_1}x_1^{\delta_2}$ for $\delta_1>0,\delta_2>0$ to get more general forms of (2.16).
\vskip.2cm Instead of a ${_pF_q}$ with argument $ax_1$ let us append $x_1^{\zeta-1}(1-x_1)^{\alpha-1}$ with a hypergeometric function ${_pF_q}$ with argument $a(1-x_1)$. Then, proceeding as before we get the normalizing constant as
$$c^{(2)}={_{p+1}F_{q+1}}(a_1,...,a_p,\alpha;b_1,...,b_q,\alpha+\zeta;a){{\Gamma(\zeta)\Gamma(\alpha)}\over{\Gamma(\alpha+\zeta)}}.\eqno(2.17)
$$Then the density of $u={{x_2}\over{x_1}}$, denoted again by $g(u)$, is given by
$$\eqalignno{g(u)&={{1}\over{c^{(2)}}}\sum_{k=0}^{\infty}{{(a_1)_k...(a_p)_k}\over{(b_1)_k...(b_q)_k}}{{a^k}\over{k!}}\cr
&\times\int_v({{v}\over{u}})^{\zeta-1}(1-{{v}\over{u}})^{\alpha+k-1}(-{{v}\over{u^2}})f(v){\rm d}v.\cr
\noalign{\hbox{The integral part is}}
u^{-\zeta-\alpha}&\int_{v=0}^uv^{\zeta}(u-v)^{\alpha-1}(1-{{v}\over{u}})^kf(v){\rm d}v.\cr
\noalign{\hbox{Hence}}
g(u)&={{u^{-\zeta-\alpha}}\over{c^{(2)}}}\int_{v=0}^uv^{\zeta}(u-v)^{\alpha-1}\cr
&\times{_pF_q}(a_1,...,a_p;b_1,...,b_q;a(1-{{v}\over{u}}))f(v){\rm d}v.&(2.18)\cr}
$$A particular case of (2.18) in terms of a ${_2F_1}$ is (2.7.3) of Mathai and Haubold (2008). This particular case was given by others earlier as a generalization of Kober operator of the first kind. The argument in the hypergeometric function could have been $a^{\delta_1}(1-x_1)^{\delta_2}$ for $\delta>0,\delta_2>0$ to produce a more general case. Also, one could have taken the argument as $a^{\delta_1}x_1^{\delta_2}(1-x_1)^{\delta_3}$. In all such cases, $g(u)$ as a density makes sense and at the same time keeping the basic structure of Kober operator of the first kind. A particular case of (2.18) is the Saigo operator of the first kind when the ${_pF_q}$ is replaced by a ${_2F_1}$. The main advantage of (2.18) is that it is a direct generalization of Kober operator of the first kind and at the same time it is a statistical density of a ratio of two statistically independently distributed random variables.

\vskip.3cm\noindent{\bf 2.3.\hskip.3cm Mellin Transform of the Generalized Kober Operator of the First Kind}

\vskip.3cm Let us consider the Mellin transform in (2.16) and (2.18). For (2.16) the basic integral to be considered is
$$\eqalignno{\int_0^{\infty}u^{s-1}&[u^{-\zeta-\alpha}\int_{v=0}^u v^{\zeta}(u-v)^{\alpha-1}{{v^k}\over{u^{k}}}{\rm d}u]{\rm d}v\cr
&=\int_{v=0}^{\infty}f(v)v^{\zeta+k}[\int_{u=v}^{\infty}u^{s-1-\zeta-\alpha-k}(u-v)^{\alpha-1}{\rm d}u]{\rm d}v\cr
&=\int_{v=0}^{\infty}f(v)v^{\zeta+k}[\int_{y=0}^{\infty}y^{\alpha-1}(y+v)^{s-1-\zeta-\alpha-k}{\rm d}y]{\rm d}v\cr
&=\int_{v=0}^{\infty}v^{s-1}f(v){\rm d}v\int_{z=0}^{\infty}z^{\alpha-1}(1+z)^{-(\alpha+\zeta+k+1-s)}{\rm d}z\cr
&=f^{*}(s){{\Gamma(\alpha)\Gamma(\zeta+1-s)}\over{\Gamma(\alpha+\zeta+1-s)}}{{(\zeta+1-s)_k}\over{(\alpha+\zeta+1-s)_k}}.\cr
\noalign{\hbox{Therefore}}
M\{g(u)\hbox{ of (2.16)};s\}&={{\Gamma(\alpha)}\over{c^{(2)}}}{{\Gamma(\zeta+1-s)}\over{\Gamma(\alpha+\zeta+1-s)}}\cr
&\times{_{p+1}F_{q+1}}a_1,...,a_p,\zeta+1-s;b_1,...,b_q,\alpha+\zeta+1-s;a)f^{*}(s)&(2.19)\cr}
$$for $\Re(\alpha)>0,\Re(s)<\Re(\zeta+1), a>0$. The basic integral to be evaluated in corresponding Mellin transform of $g(u)$ in (2.18) is the following:
$$\eqalignno{\int_{u=0}^{\infty}u^{s-1-\zeta-\alpha}&[\int_{v=0}^uv^{\zeta}(u-v)^{\alpha-1}(1-{{v}\over{u}})^kf(v){\rm d}u]{\rm d}v\cr
&=\int_{v=0}^{\infty}f(v)v^{\zeta}[\int_{u=v}^{\infty}u^{s-1-\zeta-\alpha}(u-v)^{\alpha-1}(1-{{v}\over{u}})^k{\rm d}u]{\rm d}v\cr
&=\int_{v=0}^{\infty}v^{s-1}f(v){\rm d}v\int_{z=0}^{\infty}z^{\alpha+k-1}(1+z)^{-(\alpha+\zeta+1-s+k)}{\rm d}z\cr
&=f^{*}(s){{\Gamma(\alpha)\Gamma(\zeta+1-s)}\over{\Gamma(\alpha+\zeta+1-s)}}{{(\alpha)_k}\over{(\alpha+\zeta+1-s)_k}}.\cr
\noalign{\hbox{Therefore}}
M\{g(u)\hbox{ of (2.18) };s\}&={{\Gamma(\alpha)\Gamma(\zeta+1-s)}\over{\Gamma(\alpha+\zeta+1-s)}}\cr
&\times{_{p+1}F_{q+1}}(a_1,...,a_p,\alpha;b_1,...,b_q,\alpha+\zeta+1-s;a)&(2.20)\cr}
$$for $\Re(\alpha)>0,\Re(\zeta+1-s)>0,a>0$.
\vskip.3cm\noindent{\bf 3.\hskip.3cm Riemann-Liouville Operators as Mellin Convolution}

\vskip.3cm The right sided Weyl fractional integral of order $\alpha$ is given by
$${_xW_{\infty}^{-\alpha}}f ={{1}\over{\Gamma(\alpha)}}\int_x^{\infty}(t-x)^{\alpha-1}f(t){\rm d}t,\Re(\alpha)>0.\eqno(3.1)
$$The Mellin transform is ${{\Gamma(s)}\over{\Gamma(\alpha+s)}}f^{*}(\alpha+s)$. Such a form can be generated from the Mellin convolution of a product as well as a statistical density. Consider
$$f_1(x_1)={{1}\over{\Gamma(\alpha)}}(1-x_1)^{\alpha-1},\Re(\alpha)>0.
$$Note that $f_1(x_1)$ here is a constant multiple of a statistical density. In fact, $\Gamma(\alpha+1)f_1(x_1)$ is a type-1 beta density. Let $f(x_2)$ be an arbitrary function. Let $u=x_1x_2,v=x_1,x_2={{v}\over{u}}$. Then the Mellin convolution of a product for $f(x_2)$ and $x_1^{-\alpha}f_1(x_1)$ is given by
$$\eqalignno{\int_v{{v^{-\alpha}}\over{\Gamma(\alpha)}}&(1-v)^{\alpha-1}f({{u}\over{v}}){{1}\over{v}}{\rm d}v,t={{u}\over{v}},{\rm d}v=-{{u}\over{t^2}}{\rm d}t\cr
&={{1}\over{\Gamma(\alpha)}}\int_t{{1}\over{t}}({{u}\over{t}})^{-\alpha}(1-{{u}\over{t}})^{\alpha-1}f(t){\rm d}t\cr
&={{u^{-\alpha}}\over{\Gamma(\alpha)}}\int_{t=u}^{\infty}(t-u)^{\alpha-1}f(t){\rm d}t\cr
&=u^{-\alpha}[{_xW_{\infty}^{-\alpha}}f].&(3.2)\cr}
$$
\vskip.3cm\noindent{\bf Theorem 3.1.}\hskip.3cm{\it The right sided Weyl fractional integral operator of order $\alpha$, ${_xW_{\infty}^{-\alpha}}$, can be taken as a Mellin convolution of a product of $f_1(x_1)$ and $f(x_2)$ where $f_1(x_1)$ is a constant multiple of a type-1 beta density.}

\vskip.3cm Consider the Mellin convolution of a ratio where $f_1(x_1)={{x_1^{-\alpha-1}(1-x_1)^{\alpha-1}}\over{\Gamma(\alpha)}}$ and $f_2(x_2)=x_2^{\alpha}f(x_2)$. Let $u={{x_2}\over{x_1}},v=x_2,x_1={{v}\over{u}}$. Then the Mellin convolution of a ratio is
$$\eqalignno{g(u)&=\int_vf_1({{v}\over{u}})f_2(v){{v}\over{u^2}}{\rm d}v\cr
&={{1}\over{\Gamma(\alpha)}}\int_{v=0}^u(u-v)^{\alpha-1}f(v){\rm d}v&(3.3)\cr
&={_0D_x^{-\alpha}}f\cr}
$$is the left sided Riemann-Liouville fractional integral operator of order $\alpha$. Thus, the left sided Riemann-Liouville fractional integral operator can be considered to be the Mellin convolution of a ratio. Consider the Mellin transform of $g(u)$ of (3.3). That is,
$$\eqalignno{M\{g;s\}&=\int_0^{\infty}u^{s-1}g(u){\rm d}u\cr
&=\int_{u=0}^{\infty}u^{s-1}[f_1({{v}\over{u}})f_2(v)(-{{v}\over{u^2}}){\rm d}v]{\rm d}u\cr
&=\int_{v=0}^{\infty}f_2(v)[\int_{u=v}^{\infty}u^{s-1}f_1({{v}\over{u}})(-{{v}\over{u^2}}){\rm d}u]{\rm d}v.\cr}
$$Put $x_1={{v}\over{u}},-{{v}\over{u^2}}{\rm d}u={\rm d}x_1$. Then
$$\eqalignno{\int_{u=v}^{\infty}u^{s-1}f_1({{v}\over{u}})(-{{v}\over{u^2}}){\rm d}u&=v^{s-1}\int_0^1{{1}\over{x_1^{s-1}}}f_1(x_1){\rm d}x_1\cr
&=M\{f_1;2-s\}v^{s-1}\cr}
$$for $f_1(x_1)=0$ outside the interval $[0,1]$.
$$\int_{v=0}^{\infty}v^{s-1}f_2(v){\rm d}v=M\{f_2;s\}.$$Then
$$M\{g;s\}=M\{f_1;2-s\}M\{f_2;s\}\eqno(3.4)
$$or the right side is of the form
$$\int_{x_1}\int_{x_2}({{x_2}\over{x_1}})^{s-1}f_1(x_1)f_2(x_2){\rm d}x_1\wedge{\rm d}x_2=\int_{x_1}{{1}\over{x_1^{s-1}}}f_1(x_1){\rm d}x_1\int_{x_2}x_2^{s-1}f_2(x_2){\rm d}x_2\eqno(3.5)
$$or in the form of a Mellin convolution for a ratio.

\vskip.3cm\noindent{\bf Theorem 3.2.}\hskip.3cm{\it The left sided Riemann-Liouville fractional integral is the Mellin convolution of a ratio $u={{x_2}\over{x_1}}$ when the joint function of $x_1$ and $x_2$ is of the form $f_1(x_1)f_2(x_2)$ where
$$f_1(x_1)={{x_1^{-\alpha-1}}\over{\Gamma(\alpha)}}(1-x_1)^{\alpha-1},0<x_1<1
$$and zero elsewhere, and $f_2(x_2)=x_2^{\alpha}f(x_2)$ where $f(x_2)$ is an arbitrary function, such that the Mellin transforms of $f_1(x_1)$ and $f_2(x_2)$ exist. That is,
$$g(u)=\int_vf_1({{v}\over{u}})f_2(v)(-{{v}\over{u^2}}){\rm d}v={_0D_x^{-\alpha}}f.\eqno(3.6)$$}

\vskip.3cm\noindent{\bf Note 3.1.}\hskip.3cm Note that in this case $f_1(x_1)$ is not a constant multiple of a statistical density because the exponent of $x_1$ is $-\alpha-1$ where $\Re(-\alpha)<0$. Without loss of generality, $f_2(x_2)$ can be taken as a statistical density. The Mellin transform of ${_0D_x^{-\alpha}}f$ is available in the literature, see for example Mathai and Haubold (2008).
$$M\{{_0D_x^{-\alpha}}f;s\}={{\Gamma(1-\alpha-s)}\over{\Gamma(1-s)}}f^{*}(\alpha+s), \Re(s)<1,\Re(\alpha+s)<1\eqno(3.7)
$$where $f^{*}(s)$ is the Mellin transform of $f(x)$. Thus, if $f$ is replaced by $x^{-\alpha}f$ then we have
$$M\{({_0D_x^{-\alpha}}x^{-\alpha}f)(x);s\}={{\Gamma(1-\alpha-s)}\over{\Gamma(1-s)}}f^{*}(s),\Re(s)<1,\Re(\alpha+s)<1.\eqno(3.8)
$$
\vskip.3cm\noindent{\bf Acknowledgment}

\vskip.3cm The authors would like to thank the Department of Science and Technology, Government of India, New Delhi, for the financial assistance for this work under project number SR/S4/MS:287/05.

\vskip.3cm\centerline{\bf References}

\vskip.3cm\noindent Kumar, D. (2011): P-transform, {\it Integral Transforms and Special Functions}, {\bf 22}, 603-616.

\vskip.3cm\noindent Kumar, D. and Haubold, H.J. (2010): On extended thermonuclear functions through pathway model, {\it Advances in Space Research}, {\bf 45}, 698-708.

\vskip.3cm\noindent Kumar, D. and Kilbas, A.A. (2010): Fractional calculus of P-transforms, {\it Fractional Calculus \& Applied Analysis}, {\bf 13}, 309-328.

\vskip.3cm\noindent Mathai, A.M. (2012): Generalized Kr\"{a}tzel integral and associated statistical densities, {\it International Journal of Mathematical Analysis}, {\bf 6}, 2501-2510.
\vskip.3cm\noindent Mathai, A.M. and Haubold, H.J. (2008):\hskip.3cm {\it Special Functions for Applied Scientists}, Springer, New York.
\vskip.3cm\noindent Mathai, A.M., Saxena, R.K., and Haubold, H.J. (2010): {\it The H-Function: Theory and Applications}, Springer, New York.
\vskip.3cm\noindent Nair, S.S. (2009): Pathway fractional integration operator, {\it Fractional Calculus \& Applied Analysis}, {\bf 12}, 237-252.

\bye